\documentclass{amsart}

\usepackage{amsmath}
\usepackage{amsfonts}
\usepackage{amssymb}
\usepackage{graphicx}
\usepackage{wrapfig}
\usepackage{algorithmic}
\usepackage{todonotes}
\usepackage{xcolor}
\usepackage{array}
\usepackage{extarrows}
\usepackage{wasysym}
\usepackage{dutchcal}
\usepackage{tikz}
\usepackage{pgfplots}
\usetikzlibrary{plotmarks,arrows}
\tikzset{arrow/.style={-latex, shorten >= 1ex, shorten <=1ex}}
\usepackage{pgfplotstable}
\pgfplotsset{compat=1.12}
\usetikzlibrary{external}
\usetikzlibrary{shapes.arrows}
\usepackage[algo2e,ruled, linesnumbered]{algorithm2e}
\usepackage{footmisc}
\usepackage{bbm}

\usepackage[skip=0pt, font=footnotesize, width = \linewidth]{caption}
\numberwithin{equation}{section}

\usepackage{hyperref}
\hypersetup{colorlinks=true,allcolors=blue}

\usepackage[capitalize,nameinlink]{cleveref}
\crefname{section}{section}{sections}
\crefname{subsection}{subsection}{subsections}
\Crefname{section}{Section}{Sections}
\Crefname{subsection}{Subsection}{Subsections}
\crefname{algo}{Algorithm}{Algorithms}
\crefname{table}{Table}{Tables}
\crefformat{equation}{\textup{#2(#1)#3}}
\crefrangeformat{equation}{\textup{#3(#1)#4--#5(#2)#6}}
\crefmultiformat{equation}{\textup{#2(#1)#3}}{ and \textup{#2(#1)#3}}
{, \textup{#2(#1)#3}}{, and \textup{#2(#1)#3}}
\crefrangemultiformat{equation}{\textup{#3(#1)#4--#5(#2)#6}}%
{ and \textup{#3(#1)#4--#5(#2)#6}}{, \textup{#3(#1)#4--#5(#2)#6}}{, and \textup{#3(#1)#4--#5(#2)#6}}
\Crefformat{equation}{#2Equation~\textup{(#1)}#3}
\Crefrangeformat{equation}{Equations~\textup{#3(#1)#4--#5(#2)#6}}
\Crefmultiformat{equation}{Equations~\textup{#2(#1)#3}}{ and \textup{#2(#1)#3}}
{, \textup{#2(#1)#3}}{, and \textup{#2(#1)#3}}
\Crefrangemultiformat{equation}{Equations~\textup{#3(#1)#4--#5(#2)#6}}%
{ and \textup{#3(#1)#4--#5(#2)#6}}{, \textup{#3(#1)#4--#5(#2)#6}}{, and \textup{#3(#1)#4--#5(#2)#6}}
\crefdefaultlabelformat{#2\textup{#1}#3}

\newcommand{\R}{\mathbb{R}} 
\newcommand{\N}{\mathbb{N}}

\newcommand{\tol}{\texttt{tol}}
\newcommand{\nrand}{n_\mathrm{rand}}
\newcommand{\nt}{n_\mathrm{t}}

\ifpdf
\hypersetup{
	pdftitle={DEIM vs. leverage scores for time-parallel construction of problem-adapted basis functions},
	pdfauthor={J. Schleu{\ss} and K. Smetana}
}
\fi

\title[DEIM vs. leverage scores for constructing basis functions]{DEIM vs. leverage scores for time-parallel construction of problem-adapted basis functions}

\author{Julia Schleu{\ss}}
\address{Faculty of Mathematics and Computer Science, University of M\"unster, Einsteinstr. 62, 48149 M\"unster, Germany, julia.schleuss@uni-muenster.de.}
	
\author{Kathrin Smetana}
\address{Department of Mathematical Sciences, Stevens Institute of Technology, 1 Castle Point Terrace, Hoboken, NJ 07030, United States of America, ksmetana@stevens.edu.}

\date{\today}

\thanks{The work of Julia Schleu{\ss} was funded by the Deutsche Forschungsgemeinschaft (DFG, German Research Foundation) under Germany's Excellence Strategy EXC 2044-390685587, Mathematics M\"unster: Dynamics-Geometry-Structure.}



\begin{document}
	
	\begin{abstract}
		To tackle heterogeneous time-dependent problems, an algorithm that constructs problem-adapted basis functions in an embarrassingly parallel and local manner in time has recently been proposed in [Schleuß, Smetana, ter Maat, SIAM J. Sci. Comput., 2022+]. Several simulations of the problem are performed for only few time steps in parallel by starting at different, randomly drawn start time points. For this purpose, data-dependent probability distributions that are based on the (time-dependent) data functions of the problem, such as leverage scores, are employed. In this paper, we suggest as a key new contribution to perform a deterministic time point selection based on the (discrete) empirical interpolation method (DEIM) within the proposed algorithm. In numerical experiments we investigate the performance of a DEIM based time point selection and compare it to the leverage score sampling approach.
	\end{abstract}

	\maketitle


\section{Introduction}
\label{sec:1}

Recently, problem-adapted basis functions that can be constructed in an embarrassingly parallel manner in time have been proposed to tackle heterogeneous time-dependent problems in \cite{SchSmeMaa22}. The key new idea of the approach is to select important points in time and only perform local simulations of the partial differential equation (PDE) in time in parallel by choosing these points as local end time points. The resulting basis functions are defined in space and can be combined with time-stepping schemes within multiscale or model order reduction methods. Moreover, we conjecture that they might also be relevant to construct both robust and algebraic adaptive coarse spaces within domain decomposition methods (cf. \cite{HeiSme22}).

One of the well-established tools to construct such basis functions is the proper orthogonal decomposition (POD) \cite{Sir87,BerHol93,KunVol01} that performs a singular value decomposition (SVD) on the solution of the PDE evaluated in the time grid points. However, prior to reducing, the global solution trajectory has to be computed in a sequential manner in time. As numerical experiments in \cite{SchSmeMaa22} demonstrate, the method proposed in \cite{SchSmeMaa22} can outperform the POD even in a sequential setting.

To select important points in time, data-dependent sampling strategies from randomized numerical linear algebra (NLA) \cite{DriMah16,DerMah21}, such as leverage score sampling, are employed \cite{SchSmeMaa22}. For this purpose, the time-dependent data functions are discretized and represented by a matrix, where each column of the matrix corresponds to one time point in the time grid. Usually, the aforementioned sampling methods are used to construct CUR or similar low-rank matrix decompositions (see, e.g., \cite{MahDri09}) by approximating a matrix via its columns or rows. In \cite{SorEmb16} a CUR approximate matrix decomposition based on the (discrete) empirical interpolation method (DEIM) \cite{BarMad04,ChaSor10} is proposed and numerical experiments illustrate a superior performance compared to leverage score based approaches for several test cases.

Based on these observations, the key new contribution and purpose of this manuscript is to investigate the performance of a deterministic DEIM based time point selection for the basis generation algorithm proposed in \cite{SchSmeMaa22} and compare it to the randomized time point selection via leverage score sampling employed in \cite{SchSmeMaa22}. While in \cite{SorEmb16} the authors observe that the DEIM outperforms the leverage score based approaches for the purpose of low-rank matrix approximation, we here observe that this does in general not hold true in case of time point selection.

The rest of this paper is organized as follows. First, we introduce the general problem setting in section \ref{problem_setting}. Then, we recall the parallel and local in time basis generation algorithm from \cite{SchSmeMaa22} in section \ref{parallel_in_time_algoritm} and in particular discuss the selection of relevant points in time in subsection \ref{time_point_selection}. Subsequently, we present numerical experiments comparing the DEIM and leverage score based time point selection in section \ref{numerical_experiments} and draw conclusions in section \ref{conclusions}.

\section{Problem setting}
\label{problem_setting}

We target time-dependent PDEs of the form $\partial_t u(t,x) + A(t,x) \,u(t,x) = f(t,x)$, where $A(t,x)$ denotes an elliptic differential operator that might include heterogeneous time-dependent coefficient functions. After discretizing e.g. with the finite element method in space and the implicit Euler method in time, we seek for a solution of the following discrete linear system: Given a discrete representation $\mathbf{u}_0\in \R^N$ of the initial values, the solution $\mathbf{u}_n\in \R^N$ for $n = 1,\ldots, M$ is obtained by solving
\begin{align} \label{eq_time_stepping}
(\mathbf{M}+\Delta_T \mathbf{A}_n)\,  \mathbf{u}_n = \Delta_T \mathbf{F}_n + \mathbf{M}\, \mathbf{u}_{n-1}.
\end{align}
Here, $\mathbf{M}\in \R^{N\times N}$, $\mathbf{A}_n\in \R^{N\times N}$, and $\mathbf{F}_n\in \R^{N}$ are the standard mass and stiffness matrices and right-hand side vectors and $\Delta_T$ denotes the time step size. 

We here assume that the discretization in space and time is chosen sufficiently fine such that the discretization error between the exact solution and the discrete solution of (\ref{eq_time_stepping}) is small. To ensure this, problem (\ref{eq_time_stepping}) possibly needs to be very high-dimensional and is thus computationally expensive to solve, for instance, if the coefficient functions contain fine-scale features that need to be well resolved. As a remedy, we construct suitable problem-adapted basis functions that are defined in space and combine them with the time stepping scheme.

\section{Generating problem-adapted basis functions in parallel in time}
\label{parallel_in_time_algoritm}

In the following, we briefly sketch the randomized basis generation algorithm proposed in \cite{SchSmeMaa22}. 

\begin{figure}[t]
	\centering
	\begin{tikzpicture}
	\draw[thick, ->] (0,0) -- (11.6,0);
	\draw [arrow, bend angle = 65, bend left, very thick, blue] (0.95,0) to (3.05,0);
	\draw [arrow, bend angle = 65, bend left, very thick, blue] (4.75,0) to (6.85,0);
	\draw [arrow, bend angle = 65, bend left, very thick, blue] (8.15,0) to (10.25,0);
	\draw [red, very thick] plot [only marks,mark=x, mark options={scale=2.3}] coordinates{(3,0)};
	\draw [red, very thick] plot [only marks,mark=x, mark options={scale=2.3}] coordinates{(6.8,0)};
	\draw [red, very thick] plot [only marks,mark=x, mark options={scale=2.3}] coordinates{(10.2,0)};
	\draw (0,0) node[below=2pt] {$0$};
	\draw (1,0) node[below=2pt] {$t_1$};
	\draw (3,0) node[below=2pt] {$t_2$};
	\draw (4.8,0) node[below=2pt] {$t_3$};
	\draw (6.8,0) node[below=2pt] {$t_4$};
	\draw (8.2,0) node[below=2pt] {$t_5$};
	\draw (10.2,0) node[below=2pt] {$t_6$};
	\foreach \x in {0,0.2,...,11.4} \draw (\x cm, 2.5pt) -- (\x cm, -2.5pt); 
	\draw [thick] (1,-0.5) -- (1,-0.7) -- (3,-0.7) -- (3,-0.5);
	\foreach \x in {1,1.2,...,3.0} \draw (\x cm, -0.6) -- (\x cm, -0.8); 
	\draw[very thick, blue] (2.6, -0.6) -- (2.6, -0.8); 
	\draw[very thick, blue] (2.8, -0.6) -- (2.8, -0.8);
	\draw[thick, ->] (1,-0.6) -- (1,-0.9);
	\draw[very thick, ->,blue] (3,-0.6) -- (3,-0.9);
	\draw[very thick, ->, blue] (2.4,-0.6) -- (2.4,-0.9);
	\draw (1,-1.3) node[above=-1pt] {$0$};
	\draw (2.4,-1.3) node[above=-1pt] {$k$};
	\draw (3,-1.3) node[above=-2pt] {$\nt$};
	\draw [red, very thick] plot [only marks,mark=x, mark options={scale=2.2}] coordinates{(8,-0.9)};
	\node[label=right:selected time points] at (8.1,-0.9) {};
	\end{tikzpicture}
	\caption{Sketch of Algorithm 1 in \cite{SchSmeMaa22} that constructs problem-adapted basis functions in parallel in time.}
	\label{fig_algorithm}
\end{figure}
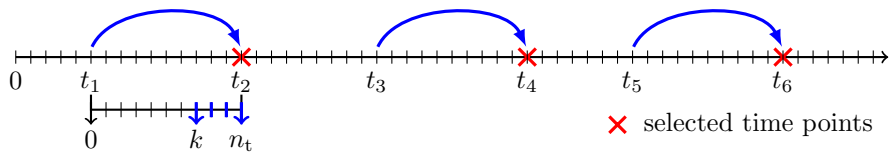

To construct suitable basis functions for the approximation of problem (\ref{eq_time_stepping}), Algorithm 1 in \cite{SchSmeMaa22} randomly chooses $\nrand \in \N$ time points in the time grid from a data-dependent probability distribution (cf. subsection \ref{time_point_selection}) and subsequently performs several simulations of the PDE for only few time steps in parallel as illustrated in Figure~\ref{fig_algorithm}. For each starting time point, a Gaussian random vector is drawn that serves as initial condition and then the local solution of the PDE is computed for $\nt$ time steps. Here, each local in time simulation of the PDE with random initial conditions approximates a local approximation space in one time point that is optimal in the sense of Kolmogorov (cf. \cite{BabLip11,BuhSme18,SchSme20,Kol36,SmePat16,MaSch21}), see section 4 and subsection 5.1 in \cite{SchSmeMaa22} for details. Subsequently, the resulting local solution trajectories evaluated at the last $\nt - k +1 $ time points are stored and compressed via an SVD to extract the reduced basis functions.

The algorithm is thus well-suited to be used on modern computer architectures as it allows to split and distribute the available computational budget over the entire time interval and facilitates an embarrassingly parallel computation.

\subsection{Discussion: Selecting points in time}
\label{time_point_selection}

In this subsection we recall how relevant points in time are randomly chosen within the basis generation algorithm in \cite{SchSmeMaa22} and propose as a key new idea a deterministic DEIM-based time point selection.

As the behavior of the solution in time is influenced by the data functions of the PDE, the time point selection in \cite{SchSmeMaa22} is guided by the behavior of the data functions in time. For this purpose, the time-dependent data functions are represented as matrices, where each column corresponds to one time point. Then, column subset selection techniques from randomized NLA \cite{DriMah16,DerMah21} are used to detect important points in time. In particular, the leverage score sampling approach \cite{DrMaMu08} appears to be an appealing choice since it is capable of detecting heterogeneous features in the data functions as numerical experiments in \cite{SchSmeMaa22} demonstrate (see subsection 5.3 in \cite{SchSmeMaa22}).

\textbf{Leverage score sampling.}\hspace{8pt}Leverage scores \cite{DrMaMu08} are computed based on the (truncated) SVD of the input matrix and reflect the influence of the individual columns or rows on the best low-rank approximation of the matrix. Therefore, leverage scores are commonly used to construct CUR or similar low-rank matrix decompositions \cite{MahDri09} that approximate a matrix directly via a subset of its columns or rows. By choosing a sufficiently large number of columns or rows the resulting approximation is nearly as good as an optimal low-rank approximation at high probability \cite{MahDri09}.

\textbf{DEIM-based selection.}\hspace{8pt}Recently, a DEIM-based CUR matrix factorization has been proposed in \cite{SorEmb16}. The method iteratively processes the leading singular vectors to deterministically choose the most important columns or rows of a matrix. In each step, the projection onto already selected components is removed from the next singular vector and a new index is selected (see Algorithm 1 in \cite{SorEmb16}). The approximation obtained by choosing exactly $r$ columns or rows is proven to be nearly as accurate as the rank-$r$ SVD \cite[Theorem 4.1]{SorEmb16}.

Numerical experiments in \cite{SorEmb16} demonstrate that the DEIM-CUR approach provides more accurate low-rank matrix approximations compared to CUR factorizations obtained from leverage scores for test cases including a sparse, non-negative matrix and matrices containing text categorization or cancer genetics data sets. The observed result holds in case of both random leverage score sampling and a deterministic selection of columns or rows with highest leverage scores. 

\textbf{New contribution.}\hspace{8pt} We here propose for the first time to employ a deterministic time point selection in Algorithm 1 of \cite{SchSmeMaa22} based on the DEIM column selection procedure introduced in \cite{SorEmb16}. We use the SVD-based DEIM approach proposed in \cite{SorEmb16} since it yields a fair comparison to the randomized leverage score sampling approach employed in \cite{SchSmeMaa22} in terms of computational costs as the (truncated) SVD is required in both approaches. Nevertheless, one could also employ a less expensive QR-based DEIM procedure as proposed in \cite{DrmGug16} and in addition further decrease computational costs by using a randomized Gram-Schmidt algorithm as introduced in \cite{BalGri22}. Alternatively, another option is to augment the DEIM approach with random sampling as proposed in \cite{PehDrmGug20}.

\section{Numerical experiments}
\label{numerical_experiments}

In the following, we numerically investigate the performance of a deterministic DEIM-based time point selection in Algorithm 1 of \cite{SchSmeMaa22} and compare it to the randomized leverage score sampling approach employed in \cite{SchSmeMaa22}. To this end, we consider two test cases from \cite{SchSmeMaa22}, i.e. problems with time-dependent source terms and a permeability coefficient that varies roughly in both space and time. These experiments give some guidance for which cases the DEIM approach might work well and for which not.

\subsection{Stove problem}

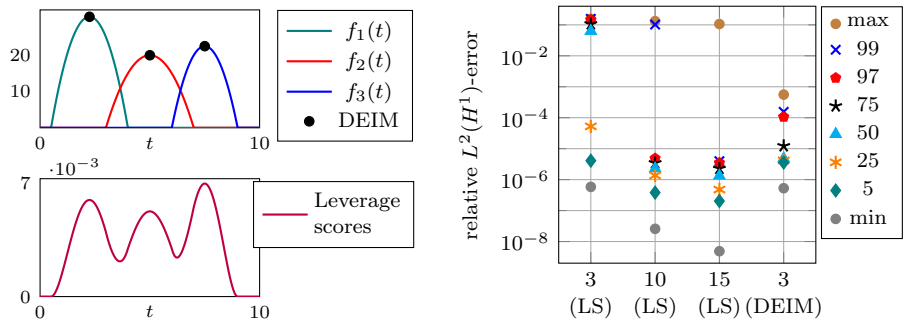
\begin{figure}
	\begin{tikzpicture}
	\begin{axis}[
	name = plot,
	title style = {font=\footnotesize, yshift = -3},
	width=4.5cm,
	height=3.15cm,
	xmin=0,
	xmax= 10,
	ymin=0,
	ymax=32.5,
	legend style={at={(1.08,1)},anchor=north west,font=\footnotesize},
	xtick = {0,5,10},
	xticklabels = {0,$t$,10},
	xtick style = {white},
	ytick={10,20},
	ytick style = {white},
	xlabel style = {yshift = 3},
	label style={font=\footnotesize},
	tick label style={font=\scriptsize}  
	]
	\addplot[solid, teal, thick] table[x index=0, y index=1] {data_figures/stove_grid_t+f_vals.dat};
	\addplot[solid, red, thick] table[x index=0, y index=2] {data_figures/stove_grid_t+f_vals.dat};
	\addplot[solid, blue, thick] table[x index=0, y index=3] {data_figures/stove_grid_t+f_vals.dat};
	\addplot[only marks, black, mark=*, mark size=1.8pt] coordinates {
		(2.25,30.6) (5,20) (7.5,22.5)
	};
	\legend{$f_1(t)$, $f_2(t)$, $f_3(t)$, DEIM}
	\end{axis}
	\begin{axis}[
	at=(plot.below south west), anchor = above north west,
	legend style={cells={align=left},at={(0.97,0.98)},anchor=north west,font=\footnotesize},
	width=4.5cm,
	height=3.15cm,
	xmin=0,
	xmax= 10,
	ymin=0,
	ymax=0.0073,
	xtick = {0,5,10},
	xticklabels = {0,$t$,10},
	xtick style = {white},
	ytick={0,0.007},
	ytick style = {white},
	label style={font=\footnotesize},
	xlabel style = {yshift = 3},
	tick label style={font=\scriptsize}  
	]
	\addplot[solid, purple, thick] table[x index=0, y index=1] {data_figures/stove_grid_t+leverage_scores.dat};
	\legend{Leverage\\scores}
	\end{axis}
	\end{tikzpicture}
	\hspace{0.25cm}
	\begin{tikzpicture}
	\begin{semilogyaxis}[
	width=5cm,
	height=5cm,
	xmin=-0.5,
	xmax=3.5,
	ymin=2e-9,
	ymax=4e-1,
	legend style={at={(1.02,1)},anchor=north west,font=\footnotesize},
	grid=both,
	grid style={line width=.1pt, draw=gray!70},
	major grid style={line width=.2pt,draw=gray!70},
	xtick={0,1,2,3},
	ytick={1e-8,1e-6, 1e-4, 1e-2},
	minor ytick={1e-7, 1e-5, 1e-3, 1e-1},
	xticklabels={3\\(LS),10\\(LS),15\\(LS), 3\\(DEIM)},
	xticklabel style = {align=center},
	ylabel= relative $L^2(H^1)$-error,
	label style={font=\footnotesize},
	tick label style={font=\footnotesize}  
	]
	\addplot[only marks, brown, mark=*, mark size=1.8pt] table[x expr = \coordindex, y index=0] {data_figures/stove_nrand_quantiles_rel_L2H1.dat};
	\addplot[only marks,blue, mark=x, mark size=2.5pt,thick] table[x expr=\coordindex, y index=1] {data_figures/stove_nrand_quantiles_rel_L2H1.dat};
	\addplot[only marks, red, mark=pentagon*, mark size=2pt] table[x expr=\coordindex, y index=3] {data_figures/stove_nrand_quantiles_rel_L2H1.dat};
	\addplot[only marks,black, mark=star,mark size=2.5pt, thick] table[x expr=\coordindex, y index=8] {data_figures/stove_nrand_quantiles_rel_L2H1.dat};
	\addplot[only marks,cyan,mark=triangle*, mark size = 2.5pt] table[x expr=\coordindex, y index=9] {data_figures/stove_nrand_quantiles_rel_L2H1.dat};
	\addplot[only marks, orange, mark=asterisk, mark size=2.5pt,thick] table[x expr=\coordindex, y index=10] {data_figures/stove_nrand_quantiles_rel_L2H1.dat};
	\addplot[only marks,teal, mark=diamond*, mark size=2.5pt] table[x expr=\coordindex, y index=11] {data_figures/stove_nrand_quantiles_rel_L2H1.dat};
	\addplot[only marks, gray, mark= oplus*, mark size=1.8pt] table[x expr=\coordindex, y index=12] {data_figures/stove_nrand_quantiles_rel_L2H1.dat};
	\legend{max,99,97,75,50,25,5,min};
	\end{semilogyaxis}
	\end{tikzpicture}
	\caption{\footnotesize Source terms $f_i$ associated with three spatially disjoint sources, DEIM time points computed from Algorithm 1 in \cite{SorEmb16} (left, top), and rank-$3$ leverage score probability distribution (left, bottom). Quantiles of relative $L^2(H^1)$-error for $\nrand=3,10,15$ time points randomly drawn from leverage scores vs. $\nrand =3$ deterministically chosen DEIM time points, $\nt = 15$, $k = 13$, $\tol = 10^{-8}$, and $100.000$ realizations of Algorithm 1 in \cite{SchSmeMaa22} (right).}
	\label{figure_stove}
\end{figure}

First, we consider the heat equation with three spatially disjoint heat sources that are turned on and off in time as illustrated in Figure \ref{figure_stove} (left, top). We refer to \cite[subsection 6.1]{SchSmeMaa22} for a full description of the numerical experiment. To choose time points within Algorithm 1 of \cite{SchSmeMaa22}, we apply the DEIM column selection Algorithm 1 from \cite{SorEmb16} to the matrix that contains the right-hand side vectors associated with each time point as columns; i.e. the same matrix that is employed for the computation of the leverage scores in \cite{SchSmeMaa22} (see Figure \ref{figure_stove} (left, bottom)).

We observe in Figure \ref{figure_stove} (left, top) that the DEIM algorithm perfectly captures the three different source terms as it selects the time points that correspond to the peaks of the stoves. Moreover, we see in Figure \ref{figure_stove} (right) that for the DEIM time point selection the relative $L^2(H^1)$-error is below $10^{-3}$ in all cases as all three stoves are detected. If we draw the same number of time points from the leverage score probability distribution, we observe that in at least $25\%$ of cases the error is below $10^{-4}$, but in at least $50\%$ of cases the error is above $5\cdot 10^{-2}$ and not all stoves are detected. For $10\, (15)$ randomly chosen time points, we see that the error is below $10^{-5}$ in $97\, (99) \%$ of cases. Hence, a small amount of oversampling is necessary to detect all three stoves at high probability when employing the randomized time point selection. Nevertheless, the corresponding local PDE simulations can be performed in an embarrassingly parallel manner and therefore do not significantly increase the wall clock time. However, we find that for this test case the deterministic DEIM approach facilitates to identify the smallest possible set of time points required for a good approximation accuracy.

\subsection{Problem with a time-dependent permeability coefficient}

\begin{figure}
	\begin{minipage}{15.5cm}
		{\footnotesize \hspace{1.3cm} $\kappa_0(x,y)$\hspace{2.2cm} solution at $t=0.1$\hspace{1.5cm} solution at $t=2$}
		\vspace{0.03cm}
		\\
		\includegraphics[width=3.95cm, height= 1cm]{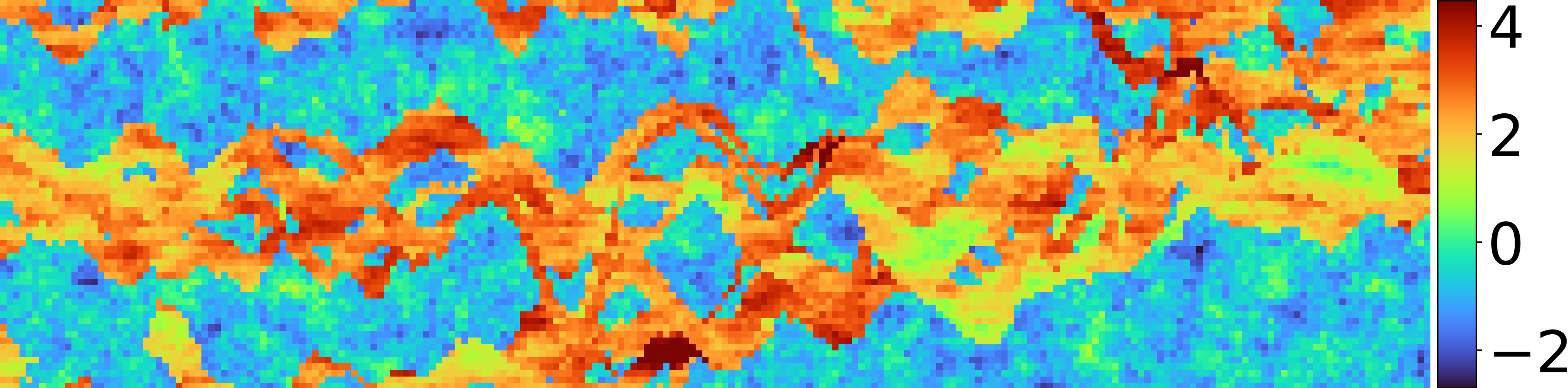}
		\includegraphics[width=3.95cm, height= 1cm]{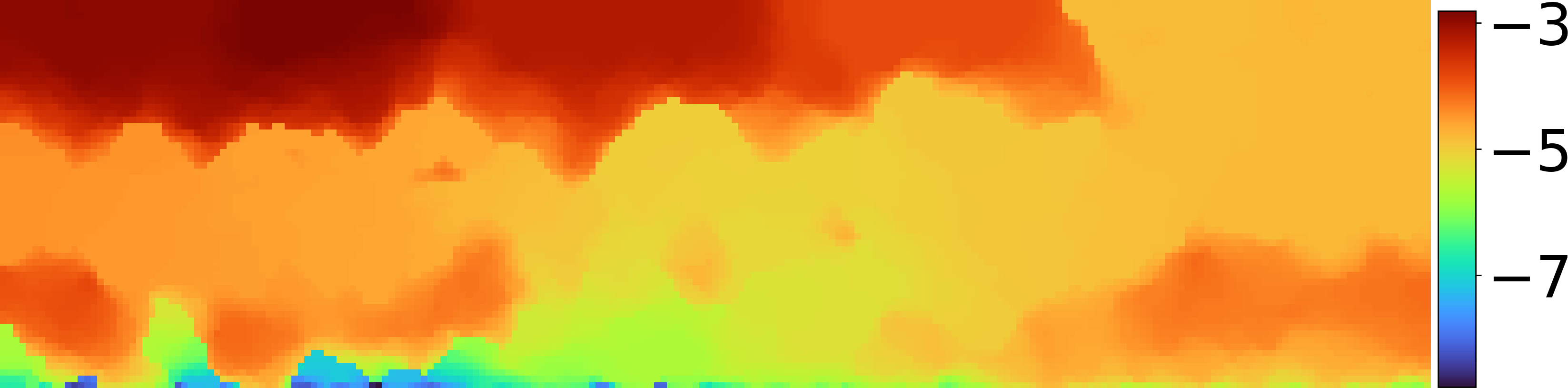}
		\includegraphics[width=3.95cm,height=1cm]{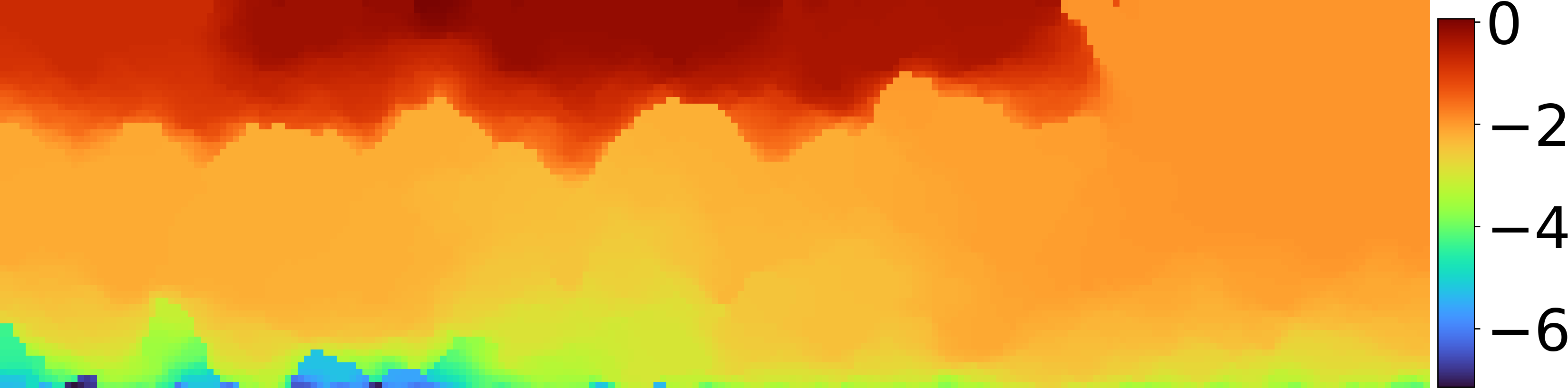}		
	\end{minipage}
	\vspace{0.15cm}\\
	\begin{minipage}{15.5cm}
		{\footnotesize \hspace{0.9cm} solution at $t=4$\hspace{1.5cm} solution at $t=7.6$\hspace{1.5cm} solution at $t=9.5$}
		\vspace{0.1cm}
		\\
		\includegraphics[width=3.95cm, height= 1cm]{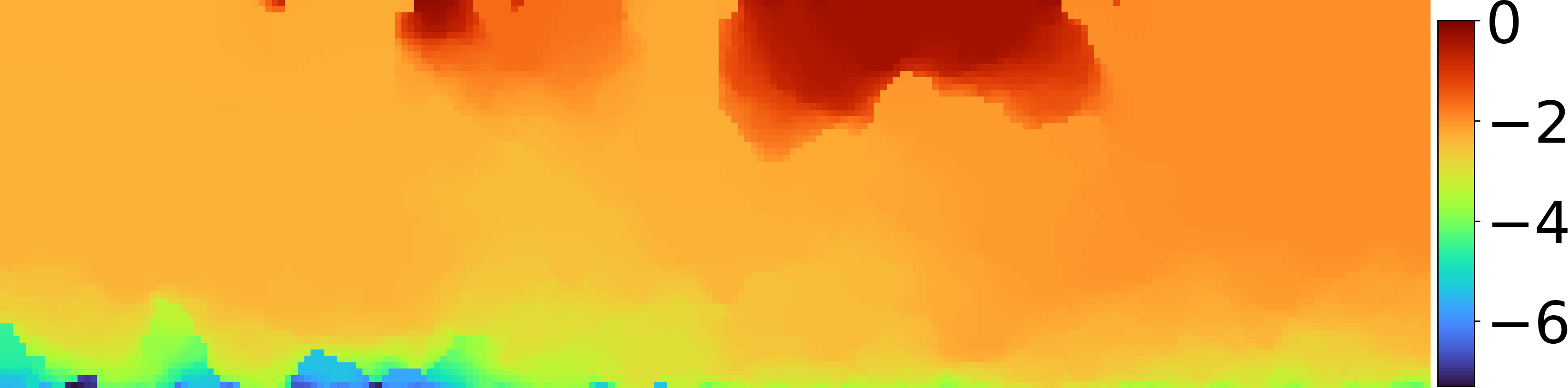}
		\includegraphics[width=3.95cm, height= 1cm]{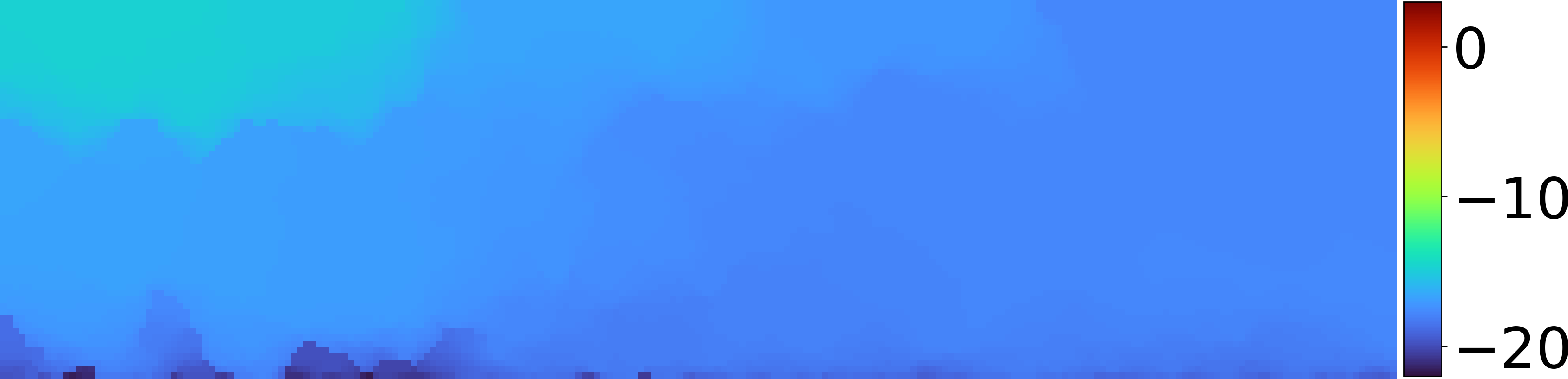}
		\includegraphics[width=3.95cm, height= 1cm]{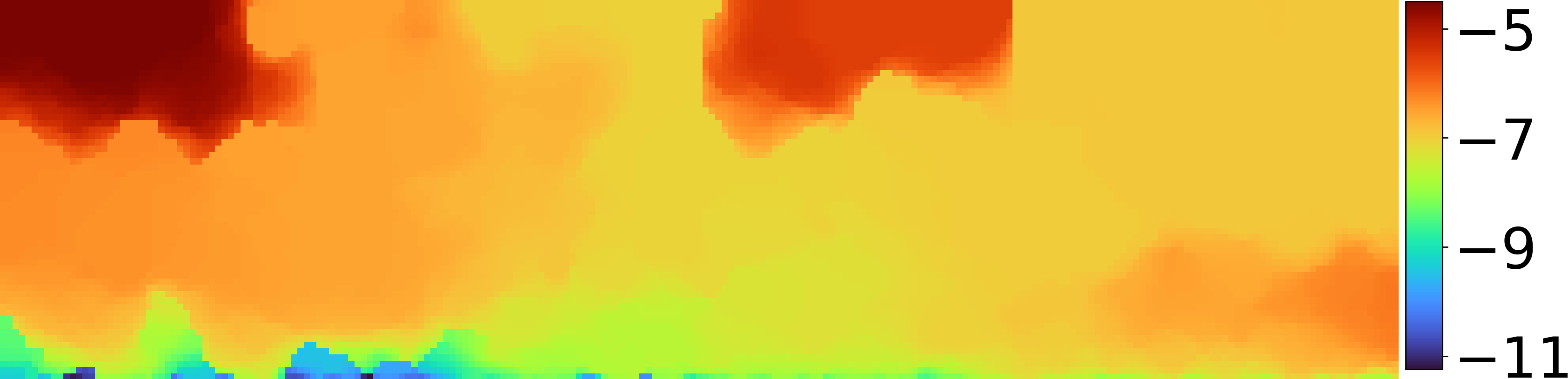}		
	\end{minipage}
	\caption{Permeability field $\kappa_0$ from \cite{SPE10} and solution evaluated at different points in time, plotted in logarithmic values to the base of $10$.}
	\label{figure_coeff_solutions_SPE10}
\end{figure}

Next, we consider an experiment including the real-world permeability coefficient $\kappa_0$ taken from the SPE10 benchmark problem \cite{SPE10}, see Figure \ref{figure_coeff_solutions_SPE10} and for details of the test case \cite[subsection 6.3]{SchSmeMaa22}. The solution trajectory of the problem is quite complex due to the different configurations and combinations of permeability and inflow into the domain depicted in Figure \ref{figure_rhs_channels_SPE10} (left and middle). As both inflow and permeability are time-dependent, we apply the DEIM column selection algorithm from \cite{SorEmb16} to both the matrix that contains the right-hand side vectors associated with each time point as columns and the matrix whose columns contain the values of the permeability field evaluated in every time point. The corresponding leverage scores are shown in Figure \ref{figure_rhs_channels_SPE10} (right).

\begin{figure}
	\begin{tikzpicture}
	\begin{axis}[
	name = MyAxis,
	width=5cm,
	height=3.2cm,
	xmin=0,
	xmax= 2.2,
	ymin=0,
	ymax=0.6,
	xtick = {1.1},
	xticklabels = $x$,
	ytick = {0.3},
	yticklabels = $y$,
	x label style={yshift=4,xshift=3},
	y label style={yshift=-12},
	title style = {font=\footnotesize, yshift = -4},
	label style={font=\footnotesize},
	tick label style={font=\scriptsize}  
	]
	\draw[fill=lightgray] (0,0) rectangle (2.2,0.6);
	\draw[fill=darkgray] (0.5,0.2) rectangle (0.6,0.6);
	\draw[fill=darkgray] (1,0.2) rectangle (1.1,0.6);
	\draw[fill=darkgray] (1.6,0.2) rectangle (1.7,0.6);
	\draw[blue, fill=blue, thick] (0.4,0.58) rectangle (1.8,0.6);
	\node[above,align = center] at (0.8,0.025) {$\kappa_1(x,y)$};
	\node[above,align = center] at (1.7,0.025) {$\kappa_2(x,y)$};
	\end{axis}
	\node[above , align=center] at (MyAxis.north) { \textcolor{blue}{$\qquad\boldsymbol{\downarrow \;\; \downarrow \;\; \downarrow \;\; \downarrow \;\; \downarrow\,\,}g_N$}};
	\end{tikzpicture}
	\hspace{0.1cm}
	\begin{tikzpicture}
	\begin{axis}[
	width=4.7cm,
	height=3.5cm,
	xmin=-0.3,
	xmax= 10.2,
	ymin=0,
	ymax=5.4,
	xtick = {},
	ytick={1,5},
	xtick = {0,5,10},
	xticklabels= {0,$t$,10},
	x label style={yshift=4},
	xtick style = {white},
	ytick style = {white},
	title style = {font=\footnotesize, yshift = -4},
	legend style={at={(0.02,1.15)},anchor=north west,font=\footnotesize},
	label style={font=\footnotesize},
	tick label style={font=\scriptsize}  
	]
	\addplot[const plot, no marks,  densely dashed, very thick, blue] coordinates {(0,0) (1,1) (5.5,0) (8,5) (9,0)} node[above,pos=.57]{};
	\addplot[const plot, no marks, thick, red] coordinates {(0,0) (3,1) (7.5,0) (8,1) (10,0)} node[above,pos=.57]{};
	\addplot[const plot, no marks, dotted, ultra thick, black] coordinates {(0,0) (8,0) (8,1) (10,1) (10,0)} node[above,pos=.57]{};
	\addplot[only marks, violet, mark=*, mark size=2pt] coordinates {
		(3,1) (8,1) (8,5) (0,0)
	};
	\legend{$g_N(t)$,$\kappa_1(t)$,$\kappa_2(t)$, DEIM},
	\end{axis}
	\end{tikzpicture}	
	\begin{tikzpicture}
	\begin{axis}[
	name=plot1,
	width = 4.7cm,
	height=3.5cm,
	xmin=-0.2,
	xmax= 10.2,
	ymin=0,
	ymax=0.0175,
	ytick={0,0.015},
	xtick = {0,5,10},
	xticklabels = {0,$t$,10},
	x label style={yshift=4},
	xtick style = {white},
	ytick style = {white},
	title style = {font=\footnotesize, yshift = -4},
	legend style={at={(0.02,0.96)},anchor=north west,font=\footnotesize},
	label style={font=\footnotesize},
	tick label style={font=\scriptsize,xshift=1.5}  
	]
	\addplot[const plot, no marks, thick, blue] coordinates {(0,0) (1,0.00066622) (5.5,0) (8,0.01665556) (9,0) (10,0)} node[above,pos=.57]{};
	\addplot[const plot, no marks, thick, purple] coordinates {(0,0) (0,0.00191571) (3,0.00147493) (7.5,0.00191571) (8,0.00330033) (10,0)} node[above,pos=.57]{};
	\legend{LS $g_N$, LS $\kappa$}
	\end{axis}
	\end{tikzpicture}
	\caption{Time-dependent inflow $g_N$ and high conductivity channels $ \kappa_1(t)\cdot \kappa_1(x,y) + \kappa_2(t)\cdot \kappa_2(x,y)$ (left and middle; dark gray equates to $10^3$, light gray to $0$). DEIM time points computed from Algorithm 1 in \cite{SorEmb16} for $g_N$ of rank-$1$ and $\kappa = \kappa_0 + \kappa_1 + \kappa_2$ of rank-$3$ (middle) and rank-$1$ leverage scores associated with $g_N$ and rank-$3$ LS corresponding to $\kappa$ (right).}
	\label{figure_rhs_channels_SPE10}
\end{figure}
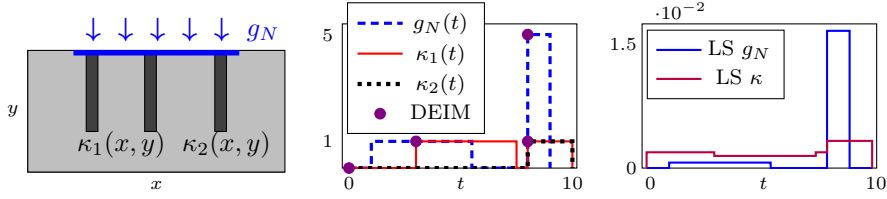

\begin{figure}
	\begin{tikzpicture}
	\begin{semilogyaxis}[
	width=4.9cm,
	height=5cm,
	xmin=-0.5,
	xmax=2.5,
	ymin=2.8e-4,
	ymax=1e0,
	legend style={at={(1.02,1)},anchor=north west,font=\footnotesize},
	grid=both,
	grid style={line width=.1pt, draw=gray!70},
	major grid style={line width=.2pt,draw=gray!70},
	label style={font=\footnotesize},
	xticklabel style ={align=center},
	xtick={0,1,2},
	ytick={1e-8,1e-7,1e-6,1e-5, 1e-4,1e-3, 1e-2,1e-1,1e0},
	minor ytick={1e-7, 1e-5, 1e-3, 1e-1},
	xticklabels={10+10\\(LS),1+3\\(DEIM), 1+3\\\;(DEIM$^\star$)},
	ylabel= relative $L^2(H^1)$-error,
	ylabel style = {yshift = -4},
	label style={font=\footnotesize},
	tick label style={font=\footnotesize} 
	]
	\addplot[only marks, brown, mark=*, mark size=1.8pt] table[x expr=\coordindex, y index=0] {data_figures/SPE10_quantiles_rel_L2H1_S.dat};
	\addplot[only marks, blue, mark=x, mark size=2.5pt,thick] table[x expr=\coordindex, y index=1] {data_figures/SPE10_quantiles_rel_L2H1_S.dat};
	\addplot[only marks, red, mark=pentagon*, mark size=2pt] table[x expr=\coordindex, y index=2] {data_figures/SPE10_quantiles_rel_L2H1_S.dat};
	\addplot[only marks,black, mark=star,mark size=2.5pt, thick] table[x expr=\coordindex, y index=3] {data_figures/SPE10_quantiles_rel_L2H1_S.dat};
	\addplot[only marks,cyan,mark=triangle*, mark size=2.5pt] table[x expr=\coordindex, y index=5] {data_figures/SPE10_quantiles_rel_L2H1_S.dat};
	\addplot[only marks, orange, mark=asterisk, mark size=2.5pt,thick] table[x expr=\coordindex, y index=7] {data_figures/SPE10_quantiles_rel_L2H1_S.dat};
	\addplot[only marks,teal, mark=diamond*, mark size=2.5pt] table[x expr=\coordindex, y index=8] {data_figures/SPE10_quantiles_rel_L2H1_S.dat};
	\addplot[only marks, gray, mark= oplus*, mark size=1.8pt] table[x expr=\coordindex, y index=9] {data_figures/SPE10_quantiles_rel_L2H1_S.dat};
	\legend{max,97,88,75,50,25,5,min};
	\end{semilogyaxis}
	\end{tikzpicture}
	\hspace{0.08cm}
	\begin{tikzpicture}
	\begin{semilogyaxis}[
	width=6cm,
	height=5.25cm,
	xmin=0,
	xmax=10,
	ymin=1e-11,
	ymax=4e0,
	legend style={at={(0.35,0.39)},anchor=north west,font=\footnotesize},
	grid=both,
	grid style={line width=.1pt, draw=gray!70},
	major grid style={line width=.2pt,draw=gray!70},
	xtick={0,2,4,6,8,10},
	xticklabels = {0,,$\qquad\quad t$,,,10},
	ytick={1e-16,1e-14,1e-12,1e-10,1e-8,1e-6, 1e-4, 1e-2,1e-0},
	minor ytick={1e-11,1e-9,1e-7, 1e-5, 1e-3, 1e-1},
	ylabel= relative $L^2(t)$-error,
	ylabel style = {yshift = -4},
	label style={font=\footnotesize},
	tick label style={font=\footnotesize}  
	]
	\addplot[densely dashdotted,purple,  very thick] table[x index = 0, y index = 1]{data_figures/SPE10_rel_L2_over_time.dat};
	\addplot[densely dotted, blue, very thick] table[x index = 0, y index = 2]{data_figures/SPE10_rel_L2_over_time.dat};
	\addplot[solid ,darkgray, thick] table[x index = 0, y index = 3]{data_figures/SPE10_rel_L2_over_time.dat};
	\legend{DEIM, DEIM$^\star$, LS}
	\end{semilogyaxis}
	\end{tikzpicture}
	\caption{Quantiles of relative $L^2(H^1)$-error for $\nrand=20 (10+10)$ time points randomly drawn from leverage scores vs. $\nrand = 4 (1+3)$ deterministically chosen DEIM time points, $\nt = 15$, $k = 13$, $\tol = 10^{-8}$, and $25.000$ realizations of Algorithm 1 in \cite{SchSmeMaa22} (left). $^\star$ indicates that time points are chosen as local start instead of end times. Mean values of relative $L^2(t)$-error (right).}
	\label{figure_quantiles_SPE10}
\end{figure}
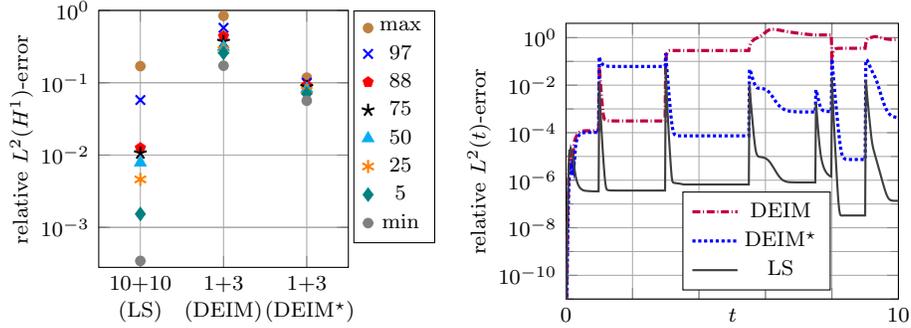

In Figure \ref{figure_quantiles_SPE10} (left) we observe that for the DEIM based time point selection the relative $L^2(H^1)$-error is above $10^{-1}$ in all cases. Moreover, the relative $L^2$-error in time is quite large in the interval $(3,10)$ as shown in Figure \ref{figure_quantiles_SPE10} (right) and we can infer that the DEIM approach is not able to detect all different configurations and combinations of the time-dependent data functions and thus all different shapes of the solution over time for this test case. This is partially due to the fact that the chosen time points are located where the data functions jump to a new value as can be seen in Figure \ref{figure_rhs_channels_SPE10} (middle). Therefore, we alternatively use the time points selected via DEIM as local start instead of end time points in Algorithm 1 in \cite{SchSmeMaa22} and observe in Figure \ref{figure_quantiles_SPE10} that both the relative $L^2(H^1)$-error and the relative $L^2$-error in the time interval $(3,10)$ decrease. However, we see that for this variant of the DEIM approach the relative $L^2(H^1)$-error is still above $5\cdot 10^{-2}$ in all cases and the solution is not captured accurately in the time interval $(1,3)$. As the rank of the time-dependent inflow (permeability) is $1$ $(3)$, we already select the maximum number of time points obtainable via the DEIM approach and the result can not be further improved unless the DEIM is, e.g., augmented with random sampling (cf. \cite{PehDrmGug20}). Moreover, we observe in Figure \ref{figure_quantiles_SPE10} (left) that the leverage score sampling approach succeeds in detecting all different shapes of the solution in at least $88\%$ of cases by drawing $20 \,(10+10)$ time points. While this is a larger number of selected time points compared to the DEIM approach, we highlight that the corresponding local PDE simulations can be performed in parallel. For this test case we thus find that the leverage score sampling approach outperforms the DEIM based time point selection.

\section{Conclusion}
\label{conclusions}

In this paper, we have proposed a deterministic DEIM based time point selection \cite{SorEmb16} within Algorithm 1 of \cite{SchSmeMaa22} that constructs problem-adapted basis functions for the approximation of (heterogeneous) time-dependent problems in an embarassingly parallel and local manner in time.

The numerical experiments demonstrate that a DEIM based time point selection can identify the smallest possible set of time points that are required to provide a good approximation accuracy. Nevertheless, we also observe that in certain cases the DEIM approach fails to detect all different combinations of the time-dependent data functions that are relevant for approximation purposes.


\begin{thebibliography}{99.}
	
	\bibitem{BabLip11}
	{ I.~Babu{\v{s}}ka and R.~Lipton}, { Optimal local approximation spaces
		for generalized finite element methods with application to multiscale
		problems}, Multiscale Model. Simul., 9 (2011), pp.~373--406.
	
	\bibitem{BalGri22}
	{ O.~Balabanov and L.~Grigori}, { Randomized {G}ram-{S}chmidt process
		with application to {GMRES}}, SIAM J. Sci. Comput., 44 (2022),
	pp.~A1450--A1474.
	
	\bibitem{BarMad04}
	{ M.~Barrault, Y.~Maday, N.~C. Nguyen, and A.~T. Patera}, { An `empirical
		interpolation' method: application to efficient reduced-basis discretization
		of partial differential equations}, C. R. Math. Acad. Sci. Paris, 339 (2004),
	pp.~667--672.
	
	\bibitem{BerHol93}
	{ G.~Berkooz, P.~Holmes, and J.~L. Lumley}, { The proper orthogonal
		decomposition in the analysis of turbulent flows}, Annu. Rev. Fluid Mech., 25
	(1993), pp.~539--575.
	
	\bibitem{BuhSme18}
	{ A.~Buhr and K.~Smetana}, { Randomized {L}ocal {M}odel {O}rder
		{R}eduction}, SIAM J. Sci. Comput., 40 (2018), pp.~A2120--A2151.
	
	\bibitem{ChaSor10}
	{ S.~Chaturantabut and D.~C. Sorensen}, { Nonlinear model reduction via
		discrete empirical interpolation}, SIAM J. Sci. Comput., 32 (2010),
	pp.~2737--2764.
	
	\bibitem{SPE10}
	{ M.~A. Christie and M.~J. Blunt}, { {Tenth SPE comparative solution
			project: A comparison of upscaling techniques}}, SPE Reservoir Evaluation \&
	Engineering, 4 (2001), pp.~308--317.
	
	\bibitem{DerMah21}
	{ M.~Derezi\'{n}ski and M.~W. Mahoney}, { Determinantal point processes
		in randomized numerical linear algebra}, Notices Amer. Math. Soc., 68 (2021),
	pp.~34--45.
	
	\bibitem{DriMah16}
	{ P.~Drineas and M.~W. Mahoney}, { {RandNLA: Randomized Numerical Linear
			Algebra}}, Commun. ACM, 59 (2016), pp.~80--90.
	
	\bibitem{DrMaMu08}
	{ P.~Drineas, M.~W. Mahoney, and S.~Muthukrishnan}, { Relative-error
		{$CUR$} matrix decompositions}, SIAM J. Matrix Anal. Appl., 30 (2008),
	pp.~844--881.
	
	\bibitem{DrmGug16}
	{ Z.~Drma\v{c} and S.~Gugercin}, { A new selection operator for the
		discrete empirical interpolation method---improved a priori error bound and
		extensions}, SIAM J. Sci. Comput., 38 (2016), pp.~A631--A648.
	
	\bibitem{HeiSme22}
	{ A.~Heinlein and K.~Smetana}, { {A fully algebraic and robust two-level
			Schwarz method based on optimal local approximation spaces}}, arXiv preprint
	arXiv:2207.05559,  (2022).
	
	\bibitem{Kol36}
	{ A.~Kolmogoroff}, { \"{U}ber die beste {A}nn{\"a}herung von {F}unktionen
		einer gegebenen {F}unktionenklasse}, Ann. of Math. (2), 37 (1936),
	pp.~107--110.
	
	\bibitem{KunVol01}
	{ K.~Kunisch and S.~Volkwein}, { Galerkin proper orthogonal decomposition
		methods for parabolic problems}, Numer. Math., 90 (2001), pp.~117--148.
	
	\bibitem{MaSch21}
	{ C.~Ma, R.~Scheichl, and T.~Dodwell}, { Novel design and analysis of
		generalized finite element methods based on locally optimal spectral
		approximations}, SIAM J. Numer. Anal., 60 (2022), pp.~244--273.
	
	\bibitem{MahDri09}
	{ M.~W. Mahoney and P.~Drineas}, { C{UR} matrix decompositions for
		improved data analysis}, Proc. Natl. Acad. Sci. U.S.A., 106 (2009),
	pp.~697--702.
	
	\bibitem{PehDrmGug20}
	{ B.~Peherstorfer, Z.~Drma\v{c}, and S.~Gugercin}, { Stability of
		discrete empirical interpolation and gappy proper orthogonal decomposition
		with randomized and deterministic sampling points}, SIAM J. Sci. Comput., 42
	(2020), pp.~A2837--A2864.
	
	\bibitem{SchSme20}
	{ J.~Schleu{\ss} and K.~Smetana}, { Optimal local approximation spaces
		for parabolic problems}, Multiscale Model. Simul., 20 (2022), pp.~551--582.
	
	\bibitem{SchSmeMaa22}
	{ J.~Schleu{\ss}, K.~Smetana, and L.~ter Maat}, { Randomized
		quasi-optimal local approximation spaces in time}, SIAM J. Sci. Comput., (2022+). (Accepted for publication. Preprint: https://arxiv.org/abs/2203.06276).
	
	\bibitem{Sir87}
	{ L.~Sirovich}, { Turbulence and the dynamics of coherent structures.
		{I}. {C}oherent structures}, Quart. Appl. Math., 45 (1987), pp.~561--571.
	
	\bibitem{SmePat16}
	{ K.~Smetana and A.~T. Patera}, { Optimal local approximation spaces for
		component-based static condensation procedures}, SIAM J. Sci. Comput., 38
	(2016), pp.~A3318--A3356.
	
	\bibitem{SorEmb16}
	{ D.~C. Sorensen and M.~Embree}, { A {DEIM} induced {CUR} factorization},
	SIAM J. Sci. Comput., 38 (2016), pp.~A1454--A1482.
	
	
	
\end{thebibliography}
\end{document}